\documentclass[12 pt]{amsart}
\usepackage{amsfonts,amsmath,amssymb,amscd}
\usepackage{xcolor}
\usepackage{comment}

\usepackage[colorlinks=true, linkcolor=blue, citecolor=blue, urlcolor=blue]{hyperref}

\newtheorem{theorem}{Theorem}[section]
\newtheorem{lemma}[theorem]{Lemma}
\newtheorem{corollary}[theorem]{Corollary}
\theoremstyle{definition}
\newtheorem{definition}[theorem]{Definition}
\newtheorem{example}[theorem]{Example}
\newtheorem{proposition}[theorem]{Proposition}

\newtheorem{qn}[theorem]{Question}

\newtheorem{remark}[theorem]{Remark}

\def \D{\mathbb{D}}
\def \C{\mathbb{C}}
\def \T{\mathbb{T}}
\def \R{\mathbb{R}}
\newcommand{\clb}{\mathcal{B}}

\newcommand{\clk}{\mathcal{K}}

\newcommand{\clf}{\mathcal{F}}

\newcommand{\na}{\mathcal{N}}
\newcommand{\ma}{\mathcal{M}}
\newcommand{\an}{\mathcal{AN}}
\newcommand{\am}{\mathcal{AM}}
\newcommand{\es}{\sigma_{ess}}
\newcommand{\vp}{\varphi}

\numberwithin{equation}{section}
\title[On $2\times 2$ block $\mathcal{AN}$ and $\mathcal{AM}$-operators]{On Absolutely norm (minimum) attaining  $2\times 2$ block operator matrix}
\author[Puspendu Nag ]{Puspendu Nag }
\address{Department of Mathematics,\\ IIT Hyderabad, Kandi,\\ Sangareddy, Telangana,\\ India, 502 284.}
\email[]{ma23resch11003@iith.ac.in}

\author[Ramesh Golla]{Ramesh Golla}
\address{Department of Mathematics\\IIT Hyderabad, Kandi\\ Sangareddy, Telangana, \\ India, 502 284.}
\email{rameshg@math.iith.ac.in}

\keywords{Norm attaining operator, Absolutely norm (minimum) attaining operator, Compact operator, Idempotents, Toeplitz operator, Hankel operator, Essential spectrum, Weyl's theorem.}

\subjclass{47A08, 47A53, 47B07, 47B35} 

\begin{document}

\begin{abstract}
    In this article, we study absolutely norm attaining operators ($\mathcal{AN}$-operators, in short), that is, operators that attain their norm on every non-zero closed subspace of a Hilbert space. Our focus is primarily on positive $2\times2$ block operator matrices in Hilbert spaces. Subsequently, we examine the analogous problem for operators that attain their minimum modulus on every nonzero closed subspace; these are referred to as absolutely minimum attaining operators (or $\mathcal{AM}$-operators, in short). We provide conditions under which these operators belong to the operator norm closure of the above two classes. In addition, we give a characterization of idempotent operators that fall into these three classes. Finally, we illustrate our results through examples that involve concrete operators.
\end{abstract}
\maketitle	

 \tableofcontents
\newpage
\section{Introduction}
This article has three main objectives. First, we investigate when a $2\times 2$ block operator matrix is absolutely norm attaining (also called an $\mathcal{AN}$-operator). Second, we study the conditions under which such a matrix is an absolutely minimum attaining operator (an $\mathcal{AM}$-operator). Third, we explore when it lies in the operator norm closure of the aforementioned classes.

A bounded linear operator $T$ on a Hilbert space  $H$  is said to be norm attaining if there exists a unit vector $x\in H_1$ such that $\|T\|=\|Tx\|$. It is said to be absolutely norm attaining if for every non-zero closed subspace $M$ of $H$, the restriction $T|_{M}:M\rightarrow H$ is norm attaining, that is, there exists a unit vector $x_M\in M$ such that $\|T|_{M}\|=\|Tx_{M}\|$. This class of operators was introduced by Carvajal and Neves in \cite{Carvajal 1}.

The class of absolutely norm attaining operators forms an important subclass of norm attaining operators, as each such operator can be expressed as a compact perturbation of a particular partial isometry (see \cite{VR} for more details). In \cite{Carvajal 1}, the authors provided several examples of $\mathcal{AN}$-operators. A complete characterization of positive $\mathcal{AN}$-operators was given in \cite{Paulsen}.  Later, characterizations for positive, self-adjoint, and normal $\mathcal{AN}$-operators were studied in \cite{VR}. A concrete representation for normal $\mathcal{AN}$-operators was presented in \cite{GRpara}, and this was subsequently generalized to hyponormal $\mathcal{AN}$-operators in \cite{NBGRhypo}. Furthermore, sufficient conditions for a hyponormal $\mathcal{AN}$-operator to be normal were discussed in the same paper. The author in \cite{Doust} gave a different proof of the characterization of positive $\mathcal{AN}$-operators and characterized weighted shifts to be in the $\mathcal{AN}$-class.  It was  shown in \cite{NBGRHISP} that a normaloid $\mathcal{AN}$-operator admits a nontrivial hyperinvariant subspace. The operator norm closure of the class of $\mathcal{AN}$-operators and several important related results were established in \cite{AN Closure}. Characterizations of Toeplitz and Hankel $\mathcal{AN}$-operators were studied in detail in \cite{CR Paris, AN Toeplitz}.

Analogous to norm attaining and absolutely norm attaining  operators, one can define minimum attaining and absolutely minimum attaining ($\mathcal{AM}$- operators) operators by replacing the norm with the minimum modulus of an operator. This class was introduced by Carvajal and Neves in \cite{Carvajal 2}. While the $\mathcal{AM}$-class bears a formal resemblance to the $\mathcal{AN}$-class, it is fundamentally different in nature. For instance, all compact operators belong to the $\mathcal{AN}$-class, whereas the $\mathcal{AM}$-class includes all finite-rank operators but excludes compact operators with infinite-dimensional range (see \cite{Carvajal 2} for more details).

The class of $\mathcal{AM}$-operators has been extensively studied. We refer the reader to the articles \cite{AM Cha, NBGRAM,CR Paris, AN Closure, NBGRHISP} for results on their spectral theory, key properties, characterizations, and various representations. It is worth noting that the class of minimum attaining operators is dense in the space of all bounded linear operators on a Hilbert space with respect to the operator norm (see \cite{SHKGRmindense} for further details).

An interesting and significant result is that the operator norm closures of the classes of absolutely norm attaining ($\mathcal{AN}$)-operators and absolutely minimum attaining ($\mathcal{AM}$)-operators coincide. This common closure has been thoroughly studied, and various properties and characterizations of operators in this class are available in the literature. Notably, every operator in this class can be expressed as a compact perturbation of a partial isometry with a finite-dimensional kernel (see \cite{AN Closure} for more details). Moreover, it has been shown (see \cite{NBGRHISP}) that any operator in this class satisfying a certain condition possesses a non-trivial hyperinvariant subspace. We refer to \cite{CR Paris, AN Toeplitz} for examples of Toeplitz and Hankel operators in this class.

We see interesting behaviors when considering block operator matrices on the direct sum of Hilbert spaces. It is possible for a block operator to fail to be $\an$ or $\am$, even when its components satisfy these properties locally. This observation motivates the development of a framework for analyzing when a positive block operator matrix is $\an$ or $\am$ based on its entries. In this article, we mainly focus on the study of various properties of $2\times 2$ positive block operator matrices of the form $\begin{bmatrix}
    A & X \\
    X^* & B
\end{bmatrix}$. We derive necessary and sufficient conditions for a special class of $2\times 2$ block operator matrices to be $\an$ and $\am$ as well. Moreover, we present a complete characterization of block operators that belong to the closure of the class $\an$. In section 6, we further study idempotent operators and describe when such operators belong to those classes. In particular, we show that if an idempotent operator $T$ is $\an$, then the \textit{Buckholtz} operator $T + T^* - I$ is also $\an$. However, the converse is not always true, in contrast to the norm-attaining case where it always holds (see \cite{NA Idempotent}). In the last section, we discuss examples involving multiplication, Toeplitz and Hankel operators in the block matrix forms.

\section{Preliminaries}
In this article, we work with bounded linear operators on infinite dimensional complex separable Hilbert spaces, which will be denoted by  $H_1, H_2$ and $H$ etc. The space of all bounded linear operators from $H_1$ into $H_2$ is denoted by $\clb (H_1, H_2)$, which forms a Banach space under the operator norm. If $H_1=H_2=H$, we write $\clb(H)$ instead of $ \clb(H_1, H_2)$. For an operator $T\in \clb (H_1, H_2)$, we denote its range and null space by $R(T)$ and $N(T)$ respectively. If $R(T)$ is a finite-dimensional subspace of $H_2$, then $T$ is said to be a \textit{finite rank} operator. The collection of all finite rank operators in $\clb(H_1, H_2)$ is denoted by $\clf(H_1, H_2)$ and particularly we write  $\clf(H)$ for $\clf(H,H)$. An operator $T\in \clb(H_1, H_2)$ is called compact if for any bounded set $B\subseteq H_1$, the closure of $T(B)\subseteq H_2$ is compact. The collection of all compact operators in $\clb(H_1, H_2)$ is denoted by $\clk(H_1, H_2)$ and particularly $\clk(H):=\clk(H, H)$. It is well known that  $\clk(H)$ is a closed ideal in $\clb(H)$ and it contains $\clf(H)$.  Moreover, $\clf(H)$ is dense in $\clk(H)$ with respect to the operator norm.

For any $T\in \clb (H_1, H_2)$, there exists unique $T^* \in \clb (H_2,H_1)$, called the adjoint of $T$, satisfying the equation 
\begin{equation*}
\langle Tx,y \rangle = \langle x, T^* y \rangle, \; \text{for every}\; x\in H_1, y\in H_2.
\end{equation*}
For $T\in\clb(H)$, we say that 
\begin{enumerate}
    \item[(i)] $T$ is self-adjoint if $T=T^*$,
    \item[(ii)] $T$ is normal if $TT^*=T^*T$,
    \item[(iii)] $T$ is idempotent if $T^2=T$,
    \item[(iv)] $T$ is orthogonal projection if $T^2=T=T^*$,
    \item[(v)] $T$ is positive if $\langle Tx, x \rangle \geq 0$ for all $x\in H$.
\end{enumerate}
We denote the class of all self-adjoint operators and positive operators in $\clb(H)$ by $\clb(H)_{sa}$ and $\clb(H)_{+}$, respectively. If $T_1, T_2 \in \clb(H)_{sa}$ then we define $T_1 \geq T_2$ iff $T_1 - T_2 \in \clb(H)_{+}$.

For $T\in \clb(H)_{+}$, there exists unique $B\in \clb(H)_{+}$ satisfying $B^2=T$, called the square root of $T$ and denoted by $B:=T^{\frac{1}{2}}$. If $T\in \clb(H_1, H_2)$, then $T^*T\in \clb(H_1)_{+}$ and the operator $|T|:=(T^*T)^{\frac{1}{2}}$ is called the modulus of $T$.

An operator $T\in\clb(H_1, H_2)$ is called \textit{norm attaining} if there exists $x_0 \in H_1$ with $\|x_0\|=1$ such that $\|Tx_0\|=\|T\|$. We denote the class of norm attaining operators by $\na(H_1, H_2)$ and particularly $\na(H):=\na(H,H)$. We have $\na(H_1, H_2)$ is dense in $\clb(H_1, H_2)$ under the operator norm (See \cite{Denseness} for the details). On the other hand, the \textit{minimum modulus} of $T$ is defined by $m(T):=\inf \{\|Tx\| : x\in H_1, \|x\|=1\}$. Note that $m(T)>0$ if and only if $R(T)$ is closed and $T$ is injective. An operator $T\in \clb(H_1, H_2)$ is said to be \textit{minimum attaining} if there exists $y_0 \in H_1$ with $\|y_0\|=1$ such that $\|Ty_0\|=m(T)$. We denote the class of minimum attaining operators in $\clb(H_1,H_2)$ by $\ma(H_1,H_2)$ and in particular $\ma(H):=\ma(H, H)$. See \cite{Carvajal 2} for further details about this class. 

Let $H= H_1 \oplus H_2 := \{(x_1, x_2) : x_1 \in H_1, x_2 \in H_2 \}$, then $H$ is called the direct sum of $H_1$ and $H_2$. Also, $H$ is a Hilbert space equipped with the inner product 
\begin{equation*}
    \langle (v_1, v_2), (w_1, w_2) \rangle_{H} = \langle v_1, w_1 \rangle_{H_1} + \langle v_2, w_2 \rangle_{H_2},
\end{equation*}
where $v_j, w_j \in H_j$ for $j=1,2.$ Let $P_j$ denote the orthogonal projection of $H$ onto $H_j$ for $j=1,2$. For any $T\in \clb (H)$, its block matrix representation with respect to the decomposition $H= H_1\oplus H_2$ is given by $T=\begin{bmatrix}
    A_{11} & A_{12} \\
    A_{21} & A_{22}
\end{bmatrix}$, where $A_{ij}\in \clb(H_j, H_i)$ for $i,j \in\{1,2\}$. The operator $T$ is positive if and only if $A_{11}\geq 0, A_{22}\geq 0$ and $A_{12}=A_{11}^{\frac{1}{2}}CA_{22}^{\frac{1}{2}}$ for some contraction $C$. In this case, $A_{21}=A^*_{12}$.  Let $A_j \in \clb(H_j)$ for $j=1,2$, then we write $A_1\oplus A_2$ to denote the diagonal block operator $\begin{bmatrix}
    A_1 & 0 \\
    0 & A_2
\end{bmatrix}$ on $H_1 \oplus H_2$. We refer to \cite{Ando,Bhat-Theertho,Osaka} for more details on block operator matrices and properties of them.

For $T\in\clb(H)$, we define the \textit{spectrum} of $T$ by $$\sigma(T):=\{\lambda \in \C : T-\lambda I \text{ is not invertible}\}.$$ It is a non-empty compact set in the complex plane. An operator $T\in\clb(H)$ is called \textit{Fredholm} if $R(T)$ is closed and $N(T)$, $N(T^*)$ are finite dimensional. The \textit{essential spectrum} of $T$ is defined by $$\es(T):=\{\mu \in \C : T-\mu I \text{ is not Fredholm}\}.$$ 

The quotient algebra $\clb(H)/\clk(H)$ is called the \textit{Calkin algebra} and we denote the canonical homomorphism from $\clb(H)$ onto $\clb(H)/\clk(H)$ by $\pi$. As a consequence of Atkinson's characterization for Fredholm operators, we have $\es(T)=\sigma(\pi(T))$. Note that $\es(T)$ is also a non-empty compact subset of $\sigma(T)$. If $T\in \clb(H)_{sa}$, by \cite[Theorem VII.11, Page 236]{Simon} we have the following:
\begin{enumerate}
    \item[(i)] $\lambda\in \es(T)$ if and only if $\lambda$ is an eigenvalue of infinite multiplicity,
    \item[(ii)] $\lambda \in \es(T)$ if and only if $\lambda$ is an accumulation point of the \textit{point spectrum} $\sigma_{p}(T):=\{\lambda\in \C : T-\lambda I \text{ is not one-one}\}$.
\end{enumerate}

If $A \in \clb(H)_{sa}$ and $K\in \clk(H)_{sa}$, then by Weyl's theorem \cite[Corollary 8.16, Page 182]{Sch}, we have
$$\es(A+K)=\es(A).$$ Also, we refer to  \cite[Theorem 2]{Kover} regarding Weyl's theorem.

\section{Block $\an$-Operators}
In this section, we discuss about characterization of $2\times 2$ block operator matrix to be absolutely norm attaining.  

\begin{definition}\cite[Definition 1.2]{Carvajal 1}
    An operator $T\in \clb(H_1, H_2)$ is called \textit{absolutely norm attaining} if it attains its norm on every non-zero closed subspace $M$ of $H_1$ i.e. $T|_{M}\in \na(M,H_2)$.  
\end{definition}
The class of absolutely norm attaining operators between $H_1$ and $H_2$ is denoted by $\an(H_1, H_2)$. We simply say an element of this class is $\an$. Also, we denote the class of positive $\an$-operators by $\an(H_1,H_2)_{+}$.  See \cite{Carvajal 1, Paulsen,VR} for examples and more details about this class. Let's recall the structure theorem for positive $\an$-operators which serves as a key tool in this section.

\begin{theorem}\cite[Theorem 5.1]{Paulsen}\label{AN Structure Paul}
     Let $H$ be a complex Hilbert space of any dimension and $T\in \clb(H)_{+}$. Then $T\in \an(H)$ if and only if $T$ has the form $\alpha I + K +F$, where $\alpha \geq 0$, $K\in \clk(H)_{+}$ and $F\in \clf(H)_{sa}$.
 \end{theorem} 

 \begin{theorem}\cite[Theorem 2.5]{VR}\label{AN Structure VR}
     Let $H$ be an infinite dimensional Hilbert space and $T\in \clb(H)$. Then the following assertions are equivalent :
     \begin{enumerate}
         \item[(i)] $T\in \an(H)_{+}$
         \item[(ii)] There exists unique  triple $(K,F,\alpha)$, where $K\in \clk(H)_{+}$, $F\in \clf(H)_{+}$ and $\alpha \geq 0$ with $0\leq F \leq \alpha I$ and $KF=0$ such that $T=K-F+\alpha I$.
     \end{enumerate}
 \end{theorem}

\begin{lemma}\cite[Corollary 2.11]{VR}\label{AN T^*T}
    Let $T\in\clb(H_1, H_2)$. Then $T\in \an(H_1, H_2)$ if and only if $T^*T\in\an(H_1)_{+}$.
\end{lemma}

\begin{lemma}\cite[Lemma 3]{Kover}\label{Kover Ess}
    Let $T\in \clb(H)$. Then $\lambda \in \es (|T|)$ if and only if $\lambda^2 \in \es(T^*T)$.
\end{lemma}

The authors feel that following two Lemmas may be available in the literature. However, they could not find a reference for the same. For the sake of completeness the proofs are provided here.

\begin{lemma}\label{Closed Range}
  Let $A\in \clb(H_1)$ and $B\in \clb(H_2)$. Then $T= A\oplus B$ has closed range if and only if each of $A$ and $B$ has closed range.
\end{lemma}

\begin{proof}
    Note that $N(T)=N(A)\oplus N(B)$. First we show that $N(T)^\perp = N(A)^\perp \oplus N(B)^\perp$. Let $(v_1, v_2)\in N(A)^\perp \oplus N(B)^\perp$. Then for any $(w_1, w_2)\in N(A) \oplus N(B)$, we get
    $$
    \langle (v_1, v_2), (w_1, w_2) \rangle = \langle v_1, w_1 \rangle_{H_1} + \langle v_2, w_2 \rangle_{H_2} = 0.
    $$
    Therefore, $(v_1, v_2) \in (N(A) \oplus N(B))^\perp = N(T)^\perp$, and hence $N(A)^\perp \oplus N(B)^\perp \subseteq N(T)^\perp$.\\

    Conversely, let $(x_1, x_2) \in N(T)^\perp$. Then for any $v_1 \in N(A)$,
     \begin{equation*}
    \langle (x_1, x_2), (v_1, 0) \rangle = 0 \implies \langle x_1, v_1 \rangle_{H_1} = 0.
    \end{equation*} 
    
    Thus, $x_1 \in N(A)^\perp$. Also, for any $v_2 \in N(B)$, 
    \begin{equation*}
        \langle (x_1, x_2), (0, v_2) \rangle = 0 \implies \langle x_2, v_2 \rangle_{H_2} = 0.
    \end{equation*}
    Thus, $x_2 \in N(B)^\perp$. Therefore, $(x_1, x_2) \in N(A)^\perp \oplus N(B)^\perp$, and hence $N(T)^\perp \subseteq N(A)^\perp \oplus N(B)^\perp$.

    Next, assume that $R(T)$ is closed. Then there exists $k > 0$ such that 
    $$ \|T(x_1, x_2) \| \geq k \|(x_1, x_2)\| \text{ for every } (x_1, x_2) \in N(T)^\perp.$$ This implies that 
    \begin{align}
    \|Ax_1 \|^2_{H_1} + \|Bx_2 \|^2_{H_2} \geq k^2 (\|x_1\|^2_{H_1} + \|x_2\|^2_{H_2})
    \end{align}
    for every  
    $x_1 \in N(A)^\perp, x_2 \in N(B)^\perp$. Choosing $x_2 =0$, we get $A$ is bounded below on $N(A)^\perp$. Again choosing $x_1=0$, we get $B$ is bounded below on $N(B)^\perp$. Thus, both $R(A)$ and $R(B)$ are closed.

    Next, assume that both $A$ and $B$ have closed range. Then there exists $k_1, k_2 > 0$ such that 
    \begin{align*}
        \|Ax_1 \|_{H_1} &\geq k_1 \|x_1 \|_{H_1} \text{ for all } x \in N(A)^\perp \\ \intertext{and }
        \|Bx_2 \|_{H_2} &\geq k_2 \|x_2 \|_{H_2} \text{ for all } x_2 \in N(B)^\perp.
    \end{align*}
    Let $\delta = \min \{k_1, k_2 \}$. Then 
        $\|Ax_1 \|^2_{H_1} + \|Bx_2 \|^2_{H_2} \geq \delta^2 (\|x_1\|^2_{H_1} + \|x_2\|^2_{H_2})$  for every $x_1 \in N(A)^\perp, x_2 \in N(B)^\perp$. Hence $\|T(x_1, x_2)\| \geq \delta \|(x_1, x_2)\|$  for every  $(x_1, x_2) \in N(A)^\perp \oplus N(B)^\perp.$
    Thus, $T$ is bounded below on $N(T)^\perp$. Hence we conclude that $T$ has closed range.
\end{proof}

\begin{lemma}\label{Diagonal Ess}
Let $A \in \clb (H_1)_{sa}$ and $B\in \clb(H_2)_{sa}$. Then $\es(A\oplus B)=\es(A)\cup \es(B)$.
\end{lemma}
\begin{proof}    
 Let $S=A\oplus B$.  Since $A=A^*$ and $B=B^*$,  $S$ is self-adjoint, hence $\sigma (S)\subseteq \R$. Then for $\lambda \in \R$, $S-\lambda I= (A-\lambda I_{H_1} )\oplus (B-\lambda I_{H_2})$.
     If $S-\lambda I$ is Fredholm, then $\dim( N(S-\lambda I))<\infty$ and $R(S-\lambda I)$ is closed. It follows that $N(A-\lambda I_{H_1})$ and $N(B-\lambda I_{H_2})$ are finite dimensional, $R(A-\lambda I_{H_1})$ and $R(B-\lambda I_{H_2})$ are closed, by Lemma \ref{Closed Range}. Thus, $\lambda\notin\es(A\oplus B)$ implies that $ \lambda \notin \es(A)$ and $\lambda \notin \es(B)$. Contrapositively, if $\lambda \in \es(A)\cup\es(B)$ then $\lambda\in\es(A\oplus B)$. Hence $\es(A)\cup\es(B)\subseteq \es(A\oplus B)$.

     Conversely, let $\lambda \in \es(A\oplus B)$. That is, either $N(S-\lambda I)$ is not finite dimensional or $R(S-\lambda I)$ is not closed.
     If $N(S-\lambda I)$ is not finite dimensional, then either  $N(A-\lambda I_{H_1})$ or $N(B-\lambda I_{H_2})$ is not finite dimensional. On the other hand, if $R(S-\lambda I)$ is not closed, then either $R(A-\lambda I_{H_1})$ or $R(B-\lambda I_{H_2})$ is not closed, by Lemma \ref{Closed Range}. Therefore, either $A-\lambda I_{H_1}$ or $B-\lambda I_{H_2}$ is not Fredholm. Thus, we get $\lambda \in \es(A) \cup\es(B)$ and hence $\es(A\oplus B) \subseteq \es(A) \cup \es(B)$.
\end{proof}

\begin{remark}
Let $A_j\in \mathcal B(H_j)$ for $j=1,2$ and $X\in \clk (H_2,H_1)$. Consider the block operator matrix $T=\begin{bmatrix}
        A_1 & X \\
        X^* & A_2
    \end{bmatrix}\in\clb (H_1 \oplus H_2)$. Then $T=\begin{bmatrix}
        A_1 & 0 \\
        0 & A_2
    \end{bmatrix} + \begin{bmatrix}
        0 & X \\
        X^* & 0
    \end{bmatrix}$ and hence by Weyl's theorem, we have $\es (T)=\es (A_1 \oplus A_2)=\es(A_1)\cup \es(A_2)$.
\end{remark}

\begin{lemma}\label{Diagonal AN}
Let $A_j\in \mathcal B(H_j)_{+}$ for $j=1,2$ and $T=\begin{bmatrix}
			A_1&0\\
			  0 & A_2
		\end{bmatrix}$. Then $T\in \an(H_1\oplus H_2)$ if and only if  $A_j \in \an(H_j)$ and $\es(A_1)=\es(A_2)$. 
\end{lemma}

\begin{proof}
    
        Since $A_j \geq 0$ for $j=1,2$, $T$ is positive.
        Let $T\in \an(H_1\oplus H_2)$. Then by Theorem \ref{AN Structure VR}, $T=\delta I + K -F$ for some $\delta \geq 0$, $K\in \clk (H_1 \oplus H_2)_{+}$ and $F\in \clf (H_1 \oplus H_2)_{+}$. Then $A_1 - \delta I_{H_1} \oplus A_2 - \delta I_{H_2} = K-F$. Let $K=\begin{bmatrix}
            K_{11} & K_{12}\\
            K_{21} & K_{22} 
        \end{bmatrix}$ and $F=\begin{bmatrix}
            F_{11} & F_{12}\\
            F_{21} & F_{22},
        \end{bmatrix}$. Since $K\geq 0,F\geq 0$, we have $K_{jj},F_{jj}\geq 0$ for $j=1,2$.
        Then 
        \begin{equation*}
        \begin{split}
            A_j - \delta I_{H_j} &=K_{jj} - F_{jj} \text{  i.e.  } A_j =\delta I_{H_j} + K_{jj} - F_{jj} \text{ for } j=1,2.
        \end{split} 
        \end{equation*}
         Hence by Theorem \ref{AN Structure Paul}, $A_j\in \an (H_j)$ with $\es (A_1)=\{\delta \}=\es (A_2)$.

         Conversely, assume that $A_j\in \an (H_j)$ for $j=1,2$. Then each of $\es(A_1),\es(A_2)$ is the singleton set, by \cite[Remark 2.6]{VR}. Let $\es(A_1)=\es(A_2)=\{\alpha\}$ for some $\alpha \geq 0$. By Lemma \ref{Diagonal Ess}, $\es (T)=\{\alpha\}$. Since $A_j\in \an (H_j)$, by Theorem \ref{AN Structure Paul}, there exists $K_j \in \clk (H_j)_{+}$ and $F_j \in \clf (H_j)_{sa}$ such that $A_j=\alpha I_{H_j} + K_j + F_j$ for $j=1,2$. Then 
         \begin{equation*}
            \begin{split}
             T&=\begin{bmatrix}
                 \alpha I_{H_1} + K_1 + F_1 & 0 \\
                 0 & \alpha I_{H_2} + K_2 + F_2
             \end{bmatrix}\\
             &=\alpha I_{H_1 \oplus H_2} + (K_1 \oplus K_2) + (F_1 \oplus F_2).
             \end{split}
         \end{equation*}
         Clearly, $K_1 \oplus K_2\in \clk (H_1 \oplus H_2)_{+}$ and $F_1 \oplus F_2 \in \clf (H_1 \oplus H_2)_{sa}$. Again by Theorem \ref{AN Structure Paul}, we get $T\in \an (H_1 \oplus H_2)$.    
\end{proof}

\begin{corollary}\label{Off Diagonal AN}
    Let $A\in \clb (H_2,H_1)$ and $B\in \clb (H_1, H_2)$. Then $T=\begin{bmatrix}
        0 & A \\
        B & 0
    \end{bmatrix} \in \an (H_1 \oplus H_2)$  if and only if $A\in \an (H_2,H_1)$, $B\in \an (H_1, H_2)$ with $\es (|A|)=\es(|B|)$.
\end{corollary}

\begin{proof}
    Note that $T^*T=\begin{bmatrix}
        0 & B^* \\
        A^* & 0
    \end{bmatrix} \begin{bmatrix}
        0 & A \\
        B & 0
    \end{bmatrix} = \begin{bmatrix}
        B^*B & 0 \\
        0 & A^*A
    \end{bmatrix}$. Now by Lemma \ref{AN T^*T}, $T\in \an (H_1 \oplus H_2)$ if and only $T^*T\in \an (H_1 \oplus H_2)_{+}$. Hence from Lemma \ref{Diagonal AN}, it follows that $B^*B \in \an (H_1)_{+}$, $A^*A \in \an (H_2)_{+}$ with $\es (B^*B)=\es (A^*A)$. Again by Lemma \ref{AN T^*T}, we get $A\in \an (H_2,H_1)$, $B\in \an (H_1, H_2)$. The condition $\es (|A|)=\es (|B|)$ follows from Lemma \ref{Kover Ess}.   
\end{proof}

\begin{theorem}\label{AN Necessary}
Let $A_j\in \mathcal B(H_j)_{+}$ for $j=1,2$ and $X\in \mathcal B(H_2,H_1)$. Let  	$T=\begin{bmatrix}
			A_1&X\\
		 X^* & A_2
		\end{bmatrix}$ be positive. If $T\in \mathcal{AN}(H_1 \oplus H_2)$, then $A_j\in \mathcal {AN}(H_j)$ with $\sigma_{ess}(A_1)=\sigma_{ess}(A_2)$ and $X\in \mathcal K(H_2,H_1)$.
        
\end{theorem}
\begin{proof}
 Assume that $T\in \an (H_1\oplus H_2)$. Then by Theorem \ref{AN Structure VR}, there exists $\alpha\geq 0, \; F\in\clf(H_1\oplus H_2)_{+}$ and $K\in \clk(H_1\oplus H_2)_{+}$ such that $T=\alpha I+K-F$.
    Suppose $K=\begin{bmatrix}
        K_{11}&K_{12}\\
        K_{21}&K_{22}
    \end{bmatrix}$ and $F=\begin{bmatrix}
        F_{11}&F_{12}\\
        F_{21}&F_{22}
    \end{bmatrix}$.
Since $K\geq 0$, by \cite[Theorem I.1.]{Ando}, we have $K_{21}=K^*_{12}$, $K_{11},K_{22}\geq 0$ and $K_{12}=K^{\frac{1}{2}}_{11}\widetilde{C}K^{\frac{1}{2}}_{22}$ for some contraction $\widetilde{C}$.
    Since $F\geq 0$, we have $F_{21}=F^*_{12}$, $F_{11},F_{22}\geq 0$ and $F_{12}=F^{\frac{1}{2}}_{11}\widehat{C}F^{\frac{1}{2}}_{22}$ for some contraction $\widehat{C}$.
    Then
    \begin{align*}
         T-\alpha I &= \begin{bmatrix}
                     K_{11}&K_{12}\\
                     K^*_{12}&K_{22}
                    \end{bmatrix} - \begin{bmatrix}
                    F_{11}&F_{12}\\
                     F^*_{12}&F_{22}
                    \end{bmatrix}\\
                    &=\begin{bmatrix}
                        K_{11}-F_{11} & K_{12}-F_{12}\\
                        K^*_{12}-F^*_{12} & K_{22}-F_{22}
                    \end{bmatrix},   
    \intertext{that is,}
        \begin{bmatrix}
            A-\alpha I_{H_1} & X\\
            X^* & B-\alpha I_{H_2}   
        \end{bmatrix}&=\begin{bmatrix}
                        K_{11}-F_{11} & K_{12}-F_{12}\\
                        K^*_{12}-F^*_{12} & K_{22}-F_{22}
                    \end{bmatrix}.                
    \end{align*}
Therefore, $A=K_{11}-F_{11}+\alpha I_{H_1}$, $B=K_{22}-F_{22}+\alpha I_{H_2}$ and $X=K_{12}-F_{12}$.
Hence, by Theorem \ref{AN Structure Paul}, $A$ and $B$ are $\an$-operators with $\es(A)=\{ \alpha \}=\es(B)$ and $X\in\mathcal K(H_2,H_1)$.
\end{proof}
Next, we discuss sufficient conditions for $T=\begin{bmatrix}
A &X\\
X^* &B
\end{bmatrix}\in \mathcal B(H_1\oplus H_2)_{+}$ to be $\mathcal{AN}$-operator. 
\begin{qn}
Let $A_j\in \mathcal{AN}(H_j)_{+}$ with $\sigma_{\text{ess}}(A_1)=\sigma_{\text{ess}}(A_2)$ for $j=1,2$. Find all compact operators $X\in \mathcal K(H_2,H_1)$ such that $T\in \mathcal{AN}(H_1\oplus H_2)_{+}$.
\end{qn}

We give a partial affirmative answer to the above question.
\begin{theorem}\label{AN Sufficient}
Let $A_j\in \mathcal B(H_j)_{+}$ for $j=1,2$ and $X\in \mathcal B(H_2,H_1)$. Let $T=\begin{bmatrix}
			A_1&X\\
		 X^* & A_2
		\end{bmatrix}$ be positive. If  $A_j\in \mathcal {AN}(H_j)$ with $\sigma_{ess}(A_1)=\sigma_{ess}(A_2)$ and $X\in \mathcal F(H_2,H_1)$, then $T\in \mathcal{AN}(H_1\oplus H_2)$.
\end{theorem}
\begin{proof}
Assume that $A_j \in \clb(H_j)$ are positive $\an$-operators with $\es(A_1)=\es(A_2)=\{\alpha\}$ and $X\in\clf(H_2,H_1)$. Then by Theorem \ref{AN Structure VR}, we have $A_j=K_{j}-F_{j}+\alpha I_{H_j}$ for some $\alpha \geq 0$ and $K_{j}\in \clk(H_1\oplus H_2)_{+}$, $F_{j}\in\clf(H_1\oplus H_2)_{+}$ for j=1,2.
Then
\begin{equation*}
    \begin{split}
        T&=\begin{bmatrix}
        A_1 & X\\
        X^* & A_2
        \end{bmatrix}\\
        &=\begin{bmatrix}
        K_1 - F_1 + \alpha I_{H_1} & X \\
        X^* & K_2 - F_2 + \alpha I_{H_2}
    \end{bmatrix}\\
        &=\begin{bmatrix}
            K_1 & 0 \\
            0 & K_2
        \end{bmatrix} + \begin{bmatrix}
            -F_1 & X\\
            X^* & -F_2
        \end{bmatrix} + \begin{bmatrix}
            \alpha I_{H_1} & 0\\
            0 & \alpha I_{H_2}
        \end{bmatrix}\\
        &=K+F+\alpha I_{H_1 \oplus H_2},
    \end{split}
\end{equation*}
where $K=\begin{bmatrix}
    K_1 & 0 \\
    0 & K_2
\end{bmatrix}$ and $F=\begin{bmatrix}
    -F_1 & X \\
    X^* & -F_2
\end{bmatrix}.$
Note that $K\geq 0$ and $F=\begin{bmatrix}
-F_1 & X\\
X^* & -F_2
\end{bmatrix}$
is self-adjoint. Hence, by Theorem \ref{AN Structure Paul},  $T\in \mathcal{AN}(H_1\oplus H_2)$.
\end{proof}

\begin{remark}
The converse of Theorem \ref{AN Sufficient} generally does not hold. That is, if $T=\begin{bmatrix}
    A & X \\
    X^* & B
\end{bmatrix}\in \an (H_1 \oplus H_2)_{+}$, then $X$ may not be always a finite rank operator. Let's illustrate with the following example.

Let $K\in \clk (H)_{+}$. Take $A=B=K+I$ and $X=-K$. Then
\begin{equation*}
    T=\begin{bmatrix}
        A & X \\
        X^* & B
    \end{bmatrix} = \begin{bmatrix}
    K+I & -K \\
    -K & K+I
\end{bmatrix}= \begin{bmatrix}
    I & 0 \\
    0 & I
\end{bmatrix} + \begin{bmatrix}
    K & - K \\
    -K & K
\end{bmatrix}.
\end{equation*}
As $-K=K^{\frac{1}{2}}(-I)K^{\frac{1}{2}}$, by \cite[Theorem I.1.]{Ando} we get $\begin{bmatrix}
    K & -K \\
    -K & K
\end{bmatrix}$ is positive, and it is also compact. Hence $T\in \an (H \oplus H)_{+}$ but here $X$ is not a finite rank operator.
\end{remark}

\begin{proposition} 
    Let $A_j\in \clb(H_j)_{+}$ for $j=1,2$ and $X\in \clb (H_2,H_1)$ with $R(X)$ closed. Also, let $T=\begin{bmatrix}
        A_1 & X \\
        X^* & A_2
    \end{bmatrix}$ be positive. Then $T\in \an (H_1 \oplus H_2)$ if and only if $A_j\in \an (H_j)$ with $\es (A_1)=\es (A_2)$ and $X\in \clf (H_2,H_1)$.
    
\end{proposition}

\begin{proof}
    Assume that $T\in \an (H_1 \oplus H_2)$. Then by Theorem \ref{AN Necessary}, we get $A_j\in \an (H_j)$ with $\es(A_1)=\es(A_2)$ and $X\in \clk (H_2,H_1)$. Since $R(X)$ is closed, it follows that $X\in \clf (H_2,H_1)$\\
    
    Conversely, let $A_j\in \an (H_j)$ with $\es(A_1)=\es(A_2)$ and $X\in \clf (H_2,H_1)$. Then by Theorem \ref{AN Sufficient}, we get $T\in \an (H_1 \oplus H_2)$.
\end{proof}

\begin{corollary}\label{General AN}
    Let $A_{ij} \in \clb (H_j,H_i)$ for $i,j=1,2$. If $A_{jj} \in \an (H_j)$ with $\es (|A_{11}|)=\es (|A_{22}|)$ and $A_{12},A_{21}$ are finite rank operators, then $T=\begin{bmatrix}
        A_{11} & A_{12} \\
        A_{21} & A_{22}
    \end{bmatrix} \in \an (H_1 \oplus H_2)$.
\end{corollary}

\begin{proof}
   Note that  
   \begin{equation*}
       \begin{split}
           T^*T&=\begin{bmatrix}
        A^*_{11} & A^*_{21} \\
        A^*_{12} & A^*_{22}
    \end{bmatrix} \begin{bmatrix}
        A_{11} & A_{12} \\
        A_{21} & A_{22}
    \end{bmatrix}\\
    &=\begin{bmatrix}
        A^*_{11}A_{11} + A^*_{21}A_{21} & A^*_{11}A_{12} + A^*_{21}A_{22} \\
        A^*_{12}A_{11}+A^*_{22}A_{21} & A^*_{12}A_{12}+A^*_{22}A_{22}
    \end{bmatrix}\\
    &=\begin{bmatrix}
        A^*_{11}A_{11} & 0 \\
        0 & A^*_{22}A_{22}
    \end{bmatrix} + \begin{bmatrix}
        A^*_{21}A_{21} & A^*_{11}A_{12} + A^*_{21}A_{22} \\
        A^*_{12}A_{11}+A^*_{22}A_{21} & A^*_{12}A_{12}
    \end{bmatrix}.
    \end{split}
    \end{equation*}

    Since $A_{jj}$ is $\an$, $A^*_{jj}A_{jj}$ is positive $\an$ for $j=1,2$, by Lemma \ref{AN T^*T}. Also, we have $\es(A^*_{11}A_{11})=\es(A^*_{22}A_{22})$ as $\es (|A_{11}|)=\es (|A_{22}|)$, by Lemma \ref{Kover Ess}. From Theorem \ref{AN Sufficient}, we get that the first block above is a positive $\an$-operator. Since $A_{12}, A_{21}$ have finite rank, the second block above is of finite rank. Thus, $T^*T$ is a positive $\an$-operator and hence by Lemma \ref{AN T^*T}, $T\in \an (H_1\oplus H_2)$.
\end{proof}

Next, we give a complete characterization of a special type of $2\times 2$ block operator matrix to be an $\mathcal{AN}$-operator.

\begin{theorem}\label{AN A_1=I}
    Let $A\in \clb (H_2)_{+}$ and $X\in \clb (H_2,H_1)$. Let $T=\begin{bmatrix}
       I_{H_1}  & X \\
        X^* & A
    \end{bmatrix} $  be positive. Then $T\in \an (H_1 \oplus H_2)$ if and only if $A\in \an (H_2)$ with $\es (A)=\{1\}$ and $X$ has finite rank.
    
\end{theorem}

\begin{proof}
   Let $T\in \an (H_1 \oplus H_2)$. Then by Theorem \ref{AN Structure VR}, there exists $\alpha \geq 0, K \in \clk (H_1 \oplus H_2)_{+}$ and $F\in \clf (H_1 \oplus H_2)_{+}$ such that
   \begin{equation*}
       T=\alpha I_{H_1 \oplus H_2} + K - F.
   \end{equation*}
    Let $K=\begin{bmatrix}
        K_{11} & K_{12} \\
        K_{21} & K_{22}
    \end{bmatrix}$ and $F=\begin{bmatrix}
        F_{11} & F_{12} \\
        F_{21} & F_{22}
    \end{bmatrix}$. Since $K\geq 0$, we have $K_{11},K_{22}\geq 0, K_{21}=K^*_{12}$ and $K_{12}=K^{\frac{1}{2}}_{11}CK^{\frac{1}{2}}_{22}$ for some contraction $C \in \clb (H_2,H_1)$. Also, $F_{21}=F^*_{12}$, since $F$ is positive. Then
    \begin{equation*}
        T=\begin{bmatrix}
            \alpha I_{H_1}+K_{11}-F_{11} & K_{12} - F_{12} \\
            K^*_{12}-F^*_{12} & \alpha I_{H_2} + K_{22} - F_{22}
        \end{bmatrix}. 
    \end{equation*}

Thus, we have  the following equations
\begin{align}
I_{H_1}=\alpha I_{H_1}+K_{11}-F_{11} \label{1}, \\
A=\alpha I_{H_2} + K_{22} - F_{22} \label{2},\\
X=K_{12}-F_{12} \label{3}.
\end{align}

Applying Weyl's theorem in \eqref{1}, we get $\alpha=1$. Consequently $K_{11}=F_{11}$ and hence $K_{12}=F^{\frac{1}{2}}_{11} C K^{\frac{1}{2}}_{22}$. Here $F_{11}$ is finite rank operator, as $F$ has finite rank. Thus, $K_{12}$ is the finite rank operator. From \eqref{3}, it follows that $X$ has finite rank. From \eqref{2} together with Theorem \ref{AN Structure Paul}, we can say that $A\in \an(H_2)$ with $\es (A)=\{\alpha \}=\{1\}$.\\

Conversely, let $A\in \an (H_2)$ with $\es (A)=\{1\}$ and $X$ has finite rank. Then we conclude that $T\in \an (H_1 \oplus H_2)$, by Theorem \ref{AN Sufficient}.  
\end{proof}

\begin{remark}\label{AN A_2=I}
    Let $A\in \clb (H_1)_{+}$ and $X\in \clb (H_2,H_1)$. Let the block operator $T=\begin{bmatrix}
        A  & X \\
        X^* & I_{H_2}
    \end{bmatrix} $ be positive. Then $T\in \an (H_1 \oplus H_2)_{+}$ if and only if $A\in \an (H_1)$ with $\es (A)=\{1\}$ and $X\in \clf (H_2,H_1)$. 
    
\end{remark}

\begin{corollary}\label{AN Contraction}

Let $C\in \mathcal B(H_2,H_1)$ be a contraction, and let $T =\begin{bmatrix}
I_{H_1}& C\\
C^* & I_{H_2}
\end{bmatrix}.$
Then $T\in \mathcal {AN}(H_1\oplus H_2)$ if and only if $C\in \mathcal F(H_2,H_1)$.

\end{corollary}

\begin{proof}
    Since $C$ is a contraction, $T$ is positive. By Theorem \ref{AN A_1=I}, it follows that $T\in \an(H_1\oplus H_2)$ if and only if $C\in\clf(H_2,H_1)$.    
\end{proof}

\section{Block $\am$-Operators}
In this section, we study when  a $2\times 2$ block operator matrix is absolutely minimum attaining. First, we recall some basic definitions and results related to this class.
\begin{definition}\cite[Definition 1.4]{Carvajal 2}
    An operator $T\in \clb(H_1, H_2)$ is said to be \textit{absolutely minimum attaining} if $T$ attains its minimum modulus on every non-zero closed subspace $M$ of $H_1$ i.e. $T\in \ma(M,H_2)$.
\end{definition}

 We denote this class by $\am(H_1, H_2)$ and call an element in this class $\am$-operator. Also, we write the class of positive $\am$-operators by $\am(H_1, H_2)_{+}$. For further details we refer to \cite{Carvajal 2,AM Cha}. We now recall the structure theorem for positive $\am$-operators that plays an essential role in this section.

\begin{theorem}\cite[Theorem 5.8]{AM Cha}\label{AM Structure}
    Let $T\in \clb(H)$. Then the following statements are equivalent :
    \begin{enumerate}
        \item[(i)] $T\in \am(H)_{+}$
        \item[(ii)] There exists $K\in\clk(H)_{+}$ with $\|K\|\leq \alpha$ and $F\in \clf(H)_{+}$ satisfying $KF=FK=0$ such that $T=\alpha I - K + F$.
    \end{enumerate}
    Also, the above decomposition is unique. 
\end{theorem}

\begin{lemma}\cite[Theorem 5.14]{AM Cha}\label{AM T^*T}
    Let $T\in \clb(H_1, H_2)$. Then the following are equivalent :
    \begin{enumerate}
    \item[(i)] $T\in\am(H_1, H_2)$
    \item[(ii)] $|T|\in\am(H_1)_{+}$
    \item[(iii)] $T^*T\in\am(H_1)_{+}$.
    \end{enumerate}
\end{lemma}

Now we show when a positive diagonal operator is $\am$ in terms of its entries.
\begin{lemma}\label{Diagonal AM}
Let $A_j\in \clb(H_j)_{+}$ for $j=1,2$ and $T=\begin{bmatrix}
			A_1&0\\
			  0 & A_2
		\end{bmatrix}$. Then $T\in \am(H_1\oplus H_2)$ if and only if  $A_j \in \am(H_j)$ and $\es(A_1)=\es(A_2)$.
\end{lemma}

\begin{proof}

Let $T\in \am (H_1 \oplus H_2)$. Since $A_j \geq 0$ for $j=1,2$, $T$ must be positive. By Theorem \ref{AM Structure}, there exists  $\alpha\geq 0$, $K\in \clk (H_1\oplus H_2)_{+}$ with $\|K\|\leq \alpha$ and $F\in \clf(H_1 \oplus H_2)_{+}$ such that $T=\alpha I_{H_1 \oplus H_2} -K+F$. Let $K=\begin{bmatrix}
    K_{11} & K_{12} \\
    K_{21} & K_{22}
\end{bmatrix}$ and $F=\begin{bmatrix}
    F_{11} & F_{12} \\
    F_{21} & F_{22}
\end{bmatrix}$.
Since $K,F \geq 0$, we have $K_{jj},F_{jj} \geq 0$ for $j=1,2$. Thus,
\begin{equation*}
    A_1=\alpha I_{H_1} - K_{11} + F_{11}, \quad A_2=\alpha I_{H_2} -K_{22} + F_{22}.
\end{equation*} 
From \cite[Proposition 3.5]{AN Toeplitz}, we get $A_{j}\in \am (H_j)$ for $j=1,2$. Also, by Weyl's theorem, $\es(A_1)=\es(A_2)=\{\alpha \}$. \\

Conversely, let $A_j\in \am(H_j)_{+}$ with $\es(A_1)=\es(A_2)$. Then by Theorem \ref{AM Structure}, there exists $K_{j}\in \clk(H_j)_{+}, F_{j}\in \clf(H_j)_{+}$ and $\alpha \geq 0$ such that $A_j=\alpha I_{H_j}-K_j+F_j$ for $j=1,2$. Then
\begin{equation*}
   \begin{split}
    T&=\begin{bmatrix}
        \alpha I_{H_1} - K_1 +F_1 & 0 \\
        0 & \alpha I_{H_2} - K_2 +F_2
    \end{bmatrix}\\
    &=\alpha I_{H_1 \oplus H_2} - (K_1 \oplus K_2)+ (F_1\oplus F_2).
    \end{split}
\end{equation*}
Clearly, $K_1 \oplus K_2 \in \clk(H_1 \oplus H_2)_{+}$ and $F_1\oplus F_2\in \clf(H_1\oplus H_2)_{+}$.  Hence by \cite[Proposition 3.5]{AN Toeplitz}, we get $T\in \am (H_1 \oplus H_2)$.
\end{proof}

By the similar steps outlined in Corollary \ref{Off Diagonal AN}, we can establish the following remark.
\begin{remark}
    Let $A\in \clb (H_2,H_1)$ and $B\in \clb (H_1, H_2)$. Then $T=\begin{bmatrix}
        0 & A \\
        B & 0
    \end{bmatrix} \in \am (H_1 \oplus H_2)$  if and only if $A\in \am (H_2,H_1)$, $B\in \am (H_1, H_2)$ with $\es (|A|)=\es(|B|)$.
\end{remark}

\begin{theorem}\label{AM Necessary}
    Let $A_j\in \clb(H_j)_{+}$ for $j=1,2$ and $X\in \clb(H_2,H_1)$. Let  	$T=\begin{bmatrix}
			A_1&X\\
		 X^* & A_2
		\end{bmatrix}\in \mathcal B(H_1\oplus H_2)$ be positive. If $T\in \am(H_1\oplus H_2)$, then $A_j\in \am(H_j)$ for $j=1,2$ with $\es(A_1)=\es(A_2)$ and $X\in \clk(H_2,H_1)$.
\end{theorem}

\begin{proof}

    Let $T\in \am (H_1 \oplus H_2)$. Then by Theorem \ref{AM Structure}, $T=\alpha I_{H_1 \oplus H_2} -K +F$ for some $\alpha \geq 0$ and $K\in \clk(H_1 \oplus H_2)_{+}, F\in \clf(H_1 \oplus H_2)_{+}$. Let $K=\begin{bmatrix}
        K_{11} & K_{12} \\
        K_{21} & K_{22}
    \end{bmatrix}$ and $F=\begin{bmatrix}
        F_{11} & F_{12} \\
        F_{21} & F_{22}
    \end{bmatrix}$. Since $K,F \geq 0$, we have $K_{jj},F_{jj}\geq 0$ for $j=1,2 ;$ $K_{21}=K^*_{12}$ and $F_{21}=F^*_{12}$. Then
    \begin{equation*}
        \begin{split}
            T&=\begin{bmatrix}
                \alpha I_{H_1} - K_{11} + F_{11} & -K_{12}+F_{12}\\
                -K^*_{12}+F^*_{12} & \alpha I_{H_2} - K_{22}+F_{22}
            \end{bmatrix}.
        \end{split}
    \end{equation*}
    This implies that  $A_j =\alpha I_{H_j} - K_{jj} + F_{jj}$ for $j=1,2$ and $X=F_{12}-K_{12}$. Hence by \cite[Proposition 3.5]{AN Toeplitz}, $A_j \in \am (H_j)$ for $j=1,2$. Also, we get $\es (A_1)=\es (A_2)=\{\alpha\}$ and $X\in \clk(H_2,H_1)$.
\end{proof}

\begin{theorem}\label{AM Sufficient}
Let $A_j\in \clb(H_j)_{+}$ for $j=1,2$ and $X\in \clb(H_2,H_1)$. Let  the block operator	$T=\begin{bmatrix}
			A_1&X\\
		 X^* & A_2
		\end{bmatrix}$ be positive on $H_1 \oplus H_2$. If  $A_j\in \am(H_j)$ with $\es(A_1)=\es(A_2)$ and $X\in \clf(H_2,H_1)$, then $T\in \am(H_1\oplus H_2)$.
\end{theorem}

\begin{proof}
Assume that $\es(A_1)=\es(A_2)=\{\alpha\}$  for some $\alpha \geq 0$. Since $A_j \in \am (H_j)$, by Theorem \ref{AM Structure}, we have $A_j = \alpha I_{H_j}-K_j+F_j$, where $\alpha \geq 0$, $K_{j} \in \clk(H_j)_{+}$ and $F_{j}\in \clf (H_j)_{+}$ for $j=1,2$.
Hence
\begin{align*}
        T&=\begin{bmatrix}
        A_1 & X\\
        X^* & A_2
        \end{bmatrix}\\
        &=\begin{bmatrix}
        \alpha I_{H_1}-K_1+F_1 & X \\
        X^* & \alpha I_{H_2}-K_2+F_2
    \end{bmatrix}\\
        &=\begin{bmatrix}
            \alpha I_{H_1} & 0 \\
            0 & \alpha I_{H_2}
        \end{bmatrix} - \begin{bmatrix}
            K_1 & 0 \\
            0 & K_2
        \end{bmatrix} + \begin{bmatrix}
            F_1 & X \\
            X^* & F_2
        \end{bmatrix}\\
        &= \alpha I_{H_1 \oplus H_2} -K_1 \oplus K_2 + \widetilde{F},
\end{align*}
where $\widetilde{F}=\begin{bmatrix}
F_1 & X\\
X^* & F_2
\end{bmatrix}$
is self-adjoint finite rank operator and $K_1 \oplus K_2$ is positive compact operator. Hence by \cite[Proposition 3.5]{AN Toeplitz}, it follows that $T\in \am (H_1 \oplus H_2)$.
\end{proof}

\begin{corollary}\label{General AM}
    Let $A_{ij} \in \clb (H_j,H_i)$ for $i,j=1,2$. If $A_{jj} \in \am (H_j)$ with $\es (|A_{11}|)=\es (|A_{22}|)$ and $A_{12},A_{21}$ are finite rank operators, then $T=\begin{bmatrix}
        A_{11} & A_{12} \\
        A_{21} & A_{22}
    \end{bmatrix} \in \am (H_1 \oplus H_2)$.
\end{corollary}

\begin{proof}
    The proof is similar to Corollary \ref{General AN}.
\end{proof}

\begin{remark}
    The converse of Theorem \ref{AM Sufficient} does not hold in general. Let $K\in \clk(H)_{+} \setminus \clf (H)_{+}$ with $\|K\|\leq \frac{1}{2}$ and assume that $A=B=I - K$ and $X=-K$. Clearly $A,B \in \am (H)_{+}$ with $\es (A)=\es (B)$. Then 
    \begin{equation*}
        T=\begin{bmatrix}
            A & X \\
            X^* & B
        \end{bmatrix} = \begin{bmatrix}
            I-K & -K \\
            -K & I-K
        \end{bmatrix} = \begin{bmatrix}
            I & 0 \\
            0 & I
        \end{bmatrix} - \begin{bmatrix}
            K & K \\
            K & K
        \end{bmatrix}.
    \end{equation*}
Since $K$ is the positive compact operator, we get $\begin{bmatrix}
    K & K \\
    K & K
\end{bmatrix} \in \clk (H\oplus H)_{+}$. Also, $-K=(I-K)^{\frac{1}{2}}C(I-K)^{\frac{1}{2}}$, where $C=-K(I-K)^{-1}$, a contraction.  Thus, by Theorem \ref{AM Structure}, $T\in \am (H\oplus H)_{+}$.
\end{remark}

\begin{corollary}
    Let $A_j\in \clb(H_j)_{+}$ for $j=1,2$ and $X\in \clb (H_2,H_1)$ with $R(X)$  closed. Let $T=\begin{bmatrix}
        A_1 & X \\
        X^* & A_2
    \end{bmatrix}$ be positive. Then $T\in \am (H_1 \oplus H_2)$ if and only if $A_j\in \am (H_j)$ with $\es (A_1)=\es (A_2)$ and $X\in \clf (H_2,H_1)$. 
    
\end{corollary}

\begin{proof}
      Let $T\in \am (H_1 \oplus H_2)$. Thus, by Theorem \ref{AM Necessary}, we get $X\in \clk(H_2,H_1)$ and $A_j\in \am(H_j)$ for $j=1,2$ with $\es(A_1)=\es(A_2)$. Since $R(X)$ is closed and $X$ is compact, so $X\in \clf (H_2,H_1)$.

The converse part clearly follows from Theorem \ref{AM Sufficient}.       
\end{proof}

Now we provide a complete characterization of a $2\times 2$ special kind of block operator to be $\am$.
\begin{theorem}\label{AM A_2=I}
     Let $A\in \clb (H_1)_{+}$ and $X\in \clb (H_2,H_1)$. Let $T=\begin{bmatrix}
        A & X \\
        X^* & I_{H_2}
    \end{bmatrix}$ be positive. Then $T\in \am(H_1\oplus H_2)$ if and only if $A\in \am (H_1)$ with $\es (A)=\{1\}$ and $X$ has finite rank. 
\end{theorem}

\begin{proof}
    Let $T\in \am (H_1\oplus H_2)$. Then by Theorem \ref{AM Necessary}, we get $A\in \am (H_1)$ with $\es(A)=\es(I_{H_2})=\{1\}$ and $X\in \clk(H_2,H_1)$. 
    
    Next, we prove that $X$ is a finite rank operator. As $T\in \am (H_1 \oplus H_2)_{+}$,  by Theorem \ref{AM Structure}, we have $T=\alpha I_{H_1\oplus H_2} - K + F$ for some $\alpha \geq 0$ and $K\in \clk(H_1\oplus H_2)_{+}$, $F\in \clf(H_1\oplus H_2)_{+}$. Also, by Weyl's theorem, $\es(T)=\{\alpha\}$. Therefore, by Lemma \ref{Diagonal Ess} we get $\alpha=1$. Let $K=\begin{bmatrix}
        K_{11} & K_{12} \\
        K_{21} & K_{22}
    \end{bmatrix}$ and $F=\begin{bmatrix}
        F_{11} & F_{12} \\
        F_{21} & F_{22}
    \end{bmatrix}$. Since $K \geq 0$, by \cite[Theorem I.1.]{Ando}, we can say $K_{21}=K^*_{12}$, $K_{jj} \geq 0$ for $j=1,2$ and $K_{12}=K_{11}^{\frac{1}{2}}CK_{22}^{\frac{1}{2}}$ for some contraction $C$. From the $(2,2)$ entry of $T$, it follows that $I_{H_2}=I_{H_2}-K_{22}+F_{22}$ that is, $K_{22}=F_{22}$. Here $F_{22}$ has finite rank since $F$ is of finite rank. Hence $K_{12}$ is a finite rank operator. Also, from $(1,2)$ entry of $T$, we get $X=F_{12}-K_{12}$, which further implies that $X\in \clf(H_2,H_1)$.
    
    The converse implication is a direct consequence of Theorem \ref{AM Sufficient}.    
\end{proof}

\begin{remark}
    Let $A\in \clb (H_2)_{+}$ and $X\in \clb (H_2,H_1)$. Assume that the block operator $T=\begin{bmatrix}
        I_{H_1} & X \\
        X^* & A
    \end{bmatrix}$ is positive. Then $T\in \am (H_1 \oplus H_2)$ if and only if $A\in \am (H_2)$ with $\es (A)=\{1\}$ and $X$ has finite rank. 
\end{remark}

\begin{corollary}
    Let $C\in \clb(H_2,H_1)$ be a contraction and $T=\begin{bmatrix}
        I_{H_1} & C \\
        C^* & I_{H_2}
    \end{bmatrix}$. Then $T\in\am (H_1 \oplus H_2)$ if and only if $C\in \clf (H_2,H_1)$.
\end{corollary}

\begin{proof}
    Since $C$ is a contraction, $T$ is a positive operator. Hence the result follows from Theorem \ref{AM A_2=I} .
\end{proof}


\section{On the closure of Block $\an $ and $\am$ Operators}

In this section, we analyze the conditions under which an operator lies in the norm-closure of the class of $\an$ as well as $\am$-operators. In addition, it is shown in \cite[Theorem 6.10]{AN Closure} that $\overline{\an(H_1, H_2)}$ and $\overline{\am(H_1, H_2)}$ coincide. Thus, without loss of generality,  it is sufficient to consider either of the class in our context.

To proceed, we recall the structure of a positive operator $T\in \clb(H)$ that lies in $\overline{\an(H)}$.
\begin{theorem}\cite[Theorem 4.2]{AN Closure}\label{AN Closure Structure}
    Let $H$ be an infinite dimensional Hilbert space and $T\in \clb(H)$. The following statements are equivalent:
    \begin{enumerate}
        \item[(i)] $T\in \overline{\an(H)}_{+}$
        \item[(ii)] There exists $\alpha \geq 0$ and $K_1,K_2\in \clk(H)_{+}$ with $K_1\leq \alpha I$ and $K_1K_2=0$ such that $T=\alpha I - K_1 + K_2$.
    \end{enumerate}
    Moreover, $K_1=(T-\alpha I)^{-}$ and $K_2=(T-\alpha I)^{+}$, indicating the decomposition of $T$ is unique. Here $A^{+}$ and $A^{-}$ denotes the positive and negative parts the self-adjoint operator $A\in \mathcal B(H)$.
\end{theorem}

\begin{lemma}\cite[Theorem 3.15]{AN Closure}\label{AN Closure T^*T}
    Let $T\in\clb(H_1, H_2)$. Then $T\in\overline{\an(H_1, H_2)}$ if and only if $T^*T\in\overline{\an(H_1)}_{+}$.
\end{lemma}

Now we establish a complete characterization of $2\times2$ block operator matrix to be in the closure of $\an$-operators.
\begin{theorem}\label{AN Closure Matrix}
Let $A_j\in \mathcal B(H_j)_{+}$ for $j=1,2$ and $X\in \mathcal B(H_2, H_1)$. Let  	$T=\begin{bmatrix}
			A_1&X\\
		 X^* & A_2
		\end{bmatrix}$ be positive. Then $T\in \overline{\mathcal{AN}(H_1\oplus H_2)}$ if and only if $A_j\in \overline{ \mathcal {AN}(H_j)}$ with $\sigma_{ess}(A_1)=\sigma_{ess}(A_2)$ and $X\in \clk(H_2, H_1)$.
\end{theorem}
\begin{proof}

Let $A_j \in \overline{\an (H_j)}$ for $j=1,2$ with $\es (A_1)=\es (A_2)$ and $X\in \clk(H_2, H_1)$. From \cite[Theorem 4.6]{AN Closure}, we have $\es (A_j)$ is the singleton, say ${\{\alpha}\}$. That is,  $\es (A_1)=\es (A_2)=\{\alpha \}$. Then by Theorem \ref{AN Closure Structure},
\begin{equation*}
    A_1=\alpha I_{H_1} - K_1 + K_2 \text{ and } A_2=\alpha I_{H_2} - K^{'}_1 + K^{'}_2,
\end{equation*}
where $K_1,K_2\in \clk(H_1)_{+}$ and $K^{'}_1,K^{'}_2 \in \clk(H_2)_{+}$ with $K_1 \leq \alpha I_{H_1}$, $K^{'}_1 \leq \alpha I_{H_2}$ and $K_1 K_2 =0=K^{'}_1 K^{'}_2$.
With this, 
\begin{equation*}
    \begin{split}
        T &=\begin{bmatrix}
            A_1 & X \\
            X^* & A_2
        \end{bmatrix}\\
        &=\begin{bmatrix}
            \alpha I_{H_1} -K_1 + K_2 & X \\
            X^* & \alpha I_{H_2} - K^{'}_1 + K^{'}_2
        \end{bmatrix}\\
        &=\begin{bmatrix}
            \alpha I_{H_1} & 0 \\
            0 & \alpha I_{H_2}
        \end{bmatrix}
        + \begin{bmatrix}
             -K_1 + K_2 & X \\ 
             X^* & - K^{'}_1 + K^{'}_2
        \end{bmatrix}\\
        &=\alpha I_{H_1 \oplus H_2} + \widetilde{K},
    \end{split}
\end{equation*}
where $\widetilde{K}=\begin{bmatrix}
             -K_1 + K_2 & X \\ 
             X^* & - K^{'}_1 + K^{'}_2
        \end{bmatrix} \in \clk (H_1 \oplus H_2)$.
From \cite[Proposition 3.5]{AN Closure}, it follows that $T\in \overline{\an (H_1 \oplus H_2)}$.

Let $T\in \overline{\an(H_1 \oplus H_2)}$. Since $T$ is positive, by Theorem \ref{AN Closure Structure}, we have 
\begin{equation*}
    T=\begin{bmatrix}
        \alpha I_{H_1} & 0 \\
        0 & \alpha I_{H_2}
    \end{bmatrix}
    - \begin{bmatrix}
        K_{11} & K_{12}\\
        K_{21} & K_{22}
    \end{bmatrix}
    +  \begin{bmatrix}
        \widehat{K}_{11} & \widehat{K}_{12}\\
        \widehat{K}_{21} & \widehat{K}_{22}
    \end{bmatrix},
\end{equation*}
where $\alpha \geq 0$, $K_{jj},\widehat{K}_{jj}\geq 0$, $K_{21}=K^*_{12}$ and $\widehat{K}_{21}=\widehat{K}^*_{12}$ .

Therefore, we get $A_j = \alpha I_{H_j} - K_{jj} + \widehat{K}_{jj}$ for $j=1,2$ and $X=\widehat{K}_{12} - K_{12} \in \clk (H_2,H_1)$. Also, by Weyl's theorem, $\es(A_1)=\es(A_2)=\{\alpha \}$, which is a singleton set. By \cite[Theorem 4.6]{AN Closure}, it follows that $A_j \in \overline{\an (H_j)}$.
\end{proof}

\begin{corollary}
Let $C\in \mathcal B(H_2,H_1)$ be a contraction, and let $T =\begin{bmatrix}
I_{H_1}& C\\
C^* & I_{H_2}
\end{bmatrix}.$
Then $T\in \overline{ \mathcal {AN}(H_1\oplus H_2)}$ if and only if $C\in \mathcal K(H_2,H_1)$.

\end{corollary}

\begin{proof}
    Since  $C$ is a contraction, we have  $T\geq 0$. Hence by Theorem \ref{AN Closure Matrix}, it follows that $T\in \overline{\an (H_1\oplus H_2)}$ if and only if $C\in \clk (H_2,H_1)$.
\end{proof}

\begin{corollary}\label{General AN Closure Matrix}
     Let $A_{ij} \in \clb (H_j,H_i)$ for $i,j=1,2$. If $A_{jj} \in \overline{\an (H_j)}$ with $\es (|A_{11}|)=\es (|A_{22}|)$ and $A_{12},A_{21}$ are compact, then $T=\begin{bmatrix}
        A_{11} & A_{12} \\
        A_{21} & A_{22}
    \end{bmatrix} \in \overline{\an (H_1 \oplus H_2)}$.
\end{corollary}

\begin{proof}
    Note that  
   \begin{equation*}
       \begin{split}
           T^*T&=\begin{bmatrix}
        A^*_{11} & A^*_{21} \\
        A^*_{12} & A^*_{22}
    \end{bmatrix} \begin{bmatrix}
        A_{11} & A_{12} \\
        A_{21} & A_{22}
    \end{bmatrix}\\
    &=\begin{bmatrix}
        A^*_{11}A_{11} + A^*_{21}A_{21} & A^*_{11}A_{12} + A^*_{21}A_{22} \\
        A^*_{12}A_{11}+A^*_{22}A_{21} & A^*_{12}A_{12}+A^*_{22}A_{22}
    \end{bmatrix}\\
    &=\widetilde{A} + \widetilde{K},
    \end{split}
    \end{equation*}
    where 
    \begin{align*}
    \widetilde{A}&=\begin{bmatrix}
        A^*_{11}A_{11} & 0 \\
        0 & A^*_{22}A_{22}
    \end{bmatrix}\\
    \intertext{and}
    \widetilde{K}&=\begin{bmatrix}
        A^*_{21}A_{21} & A^*_{11}A_{12} + A^*_{21}A_{22} \\
        A^*_{12}A_{11}+A^*_{22}A_{21} & A^*_{12}A_{12}
    \end{bmatrix}.
    \end{align*}
    Since $A_{jj}\in \overline{\an (H_j)}$, we get $A^*_{jj}A_{jj}\in \overline{\an (H_j)}_{+}$ for $j=1,2$, by Lemma \ref{AN Closure T^*T}. Also, $\es (|A_{11}|)=\es (|A_{22}|)$ implies $\es(A^*_{11} A_{11})=\es(A^*_{22} A_{22})$, by Lemma \ref{Kover Ess}. From Theorem \ref{AN Closure Matrix}, we get $\widetilde{A}\in \overline{\an (H_1 \oplus H_2)}_{+}$. Note that $\widetilde{K}$ is compact as $A_{12}$ and $A_{21}$ are compact. By \cite[Remark 4.7]{AN Closure}, it follows that $T^*T \in \overline{\an(H_1 \oplus H_2)}_{+}$ and hence $T\in \overline{\an(H_1 \oplus H_2)}$.    
\end{proof}

\section{Idempotents}

It is well known that if $P$ is an orthogonal projection on $H$, then $P\in \an (H)$ (or $P\in\am (H)$) if and only if either $N(P)$ or $R(P)$ is finite dimensional (see \cite[Theorem 3.9]{Carvajal 1}, \cite[Theorem 3.10]{Carvajal 2}). We can ask the same question for an idempotent operator. Recall that $T\in \clb (H)$ is said to be an idempotent operator if $T^2=T$. In this case, both $N(T)$ and $R(T)$ are closed. It is to be noted that $T$ is idempotent if and only if $I-T$ is idempotent. Further, we have $R(T)=N(I-T)$ and $R(I-T)=N(T)$. 

Next, we recall a geometric realization of idempotent operators. Let $T$ be an idempotent operator on $H$, and let $P$ be the orthogonal projection of $H$ onto $R(T)$. Then by \cite[Equation 1.8]{Feldman}, $T$ has the following block matrix form
$$T=\begin{bmatrix}
    I & X \\
    0 & 0
\end{bmatrix},$$ 
with respect to $H=R(P)\oplus N(P)$, where $X\in \clb(N(P),R(P))$.

\begin{theorem}\label{Idempotent AN}
Let $T$ be an idempotent operator on an infinite-dimensional Hilbert space $H$. Then the following are equivalent :
\begin{enumerate}
    \item[(i)] $T\in \an (H)$
    \item[(ii)] Either $R(T)$ or $N(T)$ is finite dimensional
    \item[(iii)] $T\in \am(H)$
    \item[(iv)] $T\in \overline{\an(H)}$.   
\end{enumerate}
\end{theorem}

\begin{proof}
$(i)\implies (ii) :$    Assume that both $N(T)$ and $R(T)$ are infinite dimensional. Then the block matrix representation of $T$ is given by
\begin{equation*}
    T=\begin{bmatrix}
        I_{R(T)} & X  \\
        0 & 0
    \end{bmatrix} \text{ on } R(T)\oplus N(T^*).
\end{equation*}
This gives us
\begin{equation*}
\begin{split}
    T^*T &=\begin{bmatrix}
        I_{R(T)} & 0 \\
        X^* & 0
    \end{bmatrix} \begin{bmatrix}
        I_{R(T)} & X  \\
        0 & 0
    \end{bmatrix}\\
    &=\begin{bmatrix}
        I_{R(T)} & X \\
        X^* & X^*X
    \end{bmatrix}.
   \end{split} 
\end{equation*}

Here $X^*X \geq 0$ on $N(T^*)$. If $T^*T \in \an (H)_{+}$, then by Theorem \ref{AN A_1=I},  $X^*X \in \an (N(T^*))$ with $\es (X^*X)=\{1\}$ and $X$ has finite rank. But $X$ being finite rank implies $X^*X$ has finite rank and so $\es (X^*X)=\{0\}$ which contradicts to $\es (X^*X)=\{1\}$. Thus, $T^*T \notin \an(H)$ and hence $T\notin \an (H)$, by Lemma \ref{AN T^*T}. Hence if $T \in \an (H)$, then either $N(T)$ or $R(T)$ is finite dimensional.

$(ii)\implies (iii) :$ Let  either $R(T)$ or $N(T)$ be finite dimensional.\\
Case I : $\dim (R(T)) < \infty$.
In this case, $T$ being a finite-rank operator is $\am$.

Case II : $\dim (N(T))<\infty$.
Since $T$ is idempotent, we have $N(T)=R(I-T)$. In this case, $I-T$ is a finite rank operator, say $F$. Then $T=I-F$, and by \cite[Lemma 3.8]{Carvajal 2}, it follows that $T\in \am(H)$.

$(iii)\implies (iv) :$ If $T\in \am(H)$, then $T\in \overline{\am(H)}=\overline{\an(H)}$, by \cite[Theorem 6.10]{AN Closure}.

$(iv)\implies (ii) :$ Let $T\in \overline{\an(H)}$ and both $R(T)$, $N(T)$ are infinite dimensional. Then by Lemma \ref{AN Closure T^*T}, $T^*T=\begin{bmatrix}
    I_{R(T)} & X \\
    X^* & X^*X
\end{bmatrix} \in \overline{\an(H)}_{+}$. From Theorem \ref{AN Closure Matrix}, we get $X\in \clk(N(T^*),R(T))$ and $\es(X^*X)=\es(I_{R(T)})=\{1\}$, which is impossible as $\es(X^*X)=\{0\}$. Thus, whenever $T\in \overline{\an(H)}$, either $R(T)$ or $N(T)$ has to be finite dimensional. 

$(ii)\implies (i) :$ If $\dim (R(T))<\infty$, then $T$ is of finite rank and so $T\in \an(H)$. On the other hand, if $\dim(N(T))<\infty$, then clearly $I-T$ is a finite rank operator. Let $I-T=F$, then $T=I-F$. Now
\begin{equation*}
\begin{split}
    T^*T &=(I-F^*)(I-F) \\
    &=I-F-F^*+F^*F\\
    &=I+F^*F - 2 Re(F).
\end{split}
\end{equation*}
Since $F^*F - 2Re(F)\in \clf (H)_{sa}$, by \cite[Proposition 4.5]{Paulsen}, it follows that $T^*T \in \an(H)$. Hence by Lemma \ref{AN T^*T}, we get $T\in \an(H)$.
\end{proof}

In \cite[Theorem 2.2]{NA Idempotent}, it is established that $T$ being idempotent on $H$ is norm attaining if and only if the \textit{Buckholtz} operator $T+T^*-I$ is  norm attaining. For the class of $\an$-operators, we show that one direction of this implication is true but its converse does not necessarily hold.

\begin{theorem}
    Let $T$ be an idempotent operator on $H$. If $T$ is $\an$, then $T+T^*-I$ is $\an$.
    \begin{proof}
        It is easy to observe that the operator $T+T^*-I$ is self-adjoint.
        Since $T=\begin{bmatrix}
            I_{R(T)} & X \\
            0 & 0
        \end{bmatrix}$ on $R(T)\oplus N(T^*)$, we get $$(T+T^*-I)^2=\begin{bmatrix}
            I_{R(T)}+XX^* & 0\\
            0 & I_{N(T^*)}+X^*X
        \end{bmatrix}.$$    
        Let $T\in \an (H)$. Then by Theorem \ref{Idempotent AN}, either $R(T)$ or $N(T)$ is finite dimensional. Also, we have $X=-P_{R(T)}|_{N(T)}(P_{N(T^*)}|_{N(T)})^{-1}$, by \cite[Theorem 2.1]{Inv idempotent}. Thus, $X$ has finite rank. Consequently $(T+T^*-I)^2$ is $\an$ as its diagonal entries are $\an$. Hence $T+T^*-I$ is the $\an$-operator. 
    \end{proof}
    
\end{theorem}

We provide a counterexample of an idempotent $T$ such that  $T+T^*-I$ is $\an$ but $T$ is not $\an$.

\begin{example}
    Consider the Hilbert space $H=l^2(\mathbb{N})$. Suppose $\{e_n : n\in \mathbb{N} \}$ is the standard orthonormal basis for $H$. Let $H_1=\bigvee \{e_{2n} : n\in \mathbb{N} \}$ and $H_2=\bigvee \{e_{2n-1} : n\in \mathbb{N} \}$. Then clearly we have $H=H_1 \oplus H_2$. Define $X: H_2\rightarrow H_1$ by  
    \begin{align*}
        X(e_{2n-1})=\begin{cases}
            \frac{1}{n} e_{2n} & \text{if } n \text{ is odd}, \\
            0 & \text{if } n \text{ is even}.
        \end{cases}
    \end{align*}
    Then $X$ is compact, both $N(X)=\bigvee\{e_{2n-1} : n \text{ is even}\}$ and $R(X)=\bigvee \{e_{2n} : n \text{ is odd}\}$ are infinite dimensional. Consider the operator $T=\begin{bmatrix}
        I_{H_1} & X \\
        0 & 0
    \end{bmatrix}$ on $H_1 \oplus H_2$. Clearly $T$ is an idempotent operator with $R(T)=H_1 \oplus \{0\}\cong H_1$ and $N(T^*)=R(T)^\perp=\{0\}\oplus H_2 \cong H_2$. Since $N(X)\subseteq H_2$ and $R(X)\subseteq H_1$, both $H_1$ and $H_2$ must be infinite dimensional. Since $T$ is idempotent, we have $N(T^*)\cong N(T)$. Therefore, we get that both $R(T)$ and $N(T)$ are infinite-dimensional. Hence by Theorem \ref{Idempotent AN}, $T\notin \an(H)$.
    
    Here $(T+T^*-I)^2=\begin{bmatrix}
        A_1 & 0 \\
        0 & A_2
    \end{bmatrix}$, where $A_1=I_{H_1} + XX^*$ and $A_2=I_{H_2} + X^*X$. Since $X$ is compact, $XX^*$ and $X^*X$ both are positive compact and hence $A_j \in \an (H_j)_{+}$ for $j=1,2$ with $\es(A_1)=\es(A_2)$. Thus, by Lemma \ref{Diagonal AN}, it follows that  $(T+T^*-I)^2\in \an(H)_{+}$, that is, $T+T^*-I \in \an (H)$.
\end{example}

\section{Examples}
In this section, we give examples to illustrate the results proved in the earlier sections.
Let $\D=\{z\in \C : |z|<1 \}$ be the open unit disk and $\T=\{z\in \C : |z|=1\}$ be the unit circle in the complex plane $\C$. Let $L^2(\T)$ be the space of square-integrable functions on $\T$ with respect to the normalized Lebesgue measure $\mu$. We denote the Hardy space by $H^2$ which is defined as $$H^2 :=\{f\in L^2(\T) : \widehat{f}(n)=0 \quad \forall n < 0 \},$$ where $\widehat{f}(n)=\langle f, \chi_{n} \rangle$ is the $n^{th}$ Fourier coefficient of $f$ and $\chi_{n} (z)=z^n$ for $z\in \T$, $n\in \mathbb{Z}$.

Let $P$ be the orthogonal projection of $L^2(\T)$ onto $H^2$. For $\vp \in L^\infty$, the \textit{multiplication operator} on $L^2(\T)$ is defined by 
\begin{equation*}
M_{\vp} f = \vp f,\; \text{ for all}\; f\in L^2(\T).
\end{equation*}
The \textit{Toeplitz operator} on $H^2$ is defined by 
\begin{equation*}
T_{\vp}f=P(\vp f),\; \text{ for all} \; f\in H^2,
\end{equation*}
and the \textit{dual Toeplitz operator} on $(H^2)^\perp :=H^2_{-}$ is defined by 
\begin{equation*}
S_{\vp}f=(I-P)f, \; \text{for all}\; f\in H^2_{-}.
\end{equation*}
The \textit{Hankel operator} $H_{\vp} : H^2 \rightarrow H^2_{-}$ is defined by 
\begin{equation*}
H_{\vp} f=(I-P)f, \; \text{ for all}\; f\in H^2.
\end{equation*}
With respect to the decomposition  $L^2(\T)=H^2 \oplus H^2_{-}$, we can express $M_\vp$ as the following block matrix form
\begin{equation*}
    M_{\vp}=\begin{bmatrix}
        T_{\vp} & H^*_{\bar{\vp}} \\
        H_{\vp} & S_{\vp}
    \end{bmatrix}.
\end{equation*}

\begin{definition}
    An operator $V : H_1 \rightarrow H_2$ is called \textit{anti-unitary} if for every $x\in H_1, y\in H_2$ and $\alpha,\beta \in \C$
    \begin{enumerate}
        \item[(i)] $\langle Vx,Vy \rangle = \langle y,x \rangle$,
        \item[(ii)] $V(\alpha x + \beta y)=\bar{\alpha}V(x) + \bar{\beta}V(y)$.
    \end{enumerate}       
\end{definition}
Let $A_{j}\in \clb(H_j)$ for $j=1,2$. We say $A_1$ and $A_2$ are anti-unitarily equivalent if there exists an anti-unitary operator $V\in \clb (H_1, H_2)$ such that $VA_1=A_2 V$.

Define $V: L^2(\T)\rightarrow L^2(\T)$ by

\begin{equation*}
Vf(z)=\bar{z}\overline{f(z)},\;  \text{for all}\; f\in L^2(\T),z\in \T.
\end{equation*}
Then $V$ is an anti-unitary operator on $L^2(\T)$. In fact, we have $V=V^*=V^{-1}$ and $VT_{\vp}=S_{\bar{\vp}}V$ i.e. $T_{\vp}$ and $S_{\bar{\vp}}$ are anti-unitarily equivalent. We refer to \cite{Dual TO} for further details.

\begin{remark}\label{Ess TO and DTO}
From the anti-unitary equivalence between $T_\vp$ and $S_{\bar{\vp}}$, we get the follwing
\begin{equation*}
   \begin{split}
    \lambda \in \es(T_\vp) &\iff T_\vp -\lambda I \text{ is not Fredholm}\\
    &\iff VS_{\bar{\vp}}V - V\bar{\lambda}V \text{ is not Fredholm}\\
    &\iff V(S^*_\vp - \bar{\lambda} I)V \text{ is not Fredholm}\\
    &\iff S^*_\vp - \bar{\lambda} I \text{ is not Fredholm}\\
    &\iff S_\vp - \lambda I \text{ is not Fredholm}\\
    &\iff \lambda \in \es(S_\vp).
    \end{split}    
\end{equation*}
Apparently we can say $\es(T_\vp)=\es(S_\vp)$. In a similar way, it can be shown that $\es(T_\vp T^*_\vp)=\es(S^*_\vp S_\vp)$ and consequently $\es(|T_{\bar{\vp}}|)=\es(|S_\vp|)$, by Lemma \ref{Kover Ess}.
\end{remark}

\begin{remark}\label{Anti}
Since $V$ is an isometry, from \cite[Proposition 3.2]{Carvajal 1} and \cite[Proposition 3.3]{Carvajal 2}, we can prove the following :
\begin{enumerate}
\item[(i)] $T_\vp$ is $\an$ if and only if $S_{\bar{\vp}}$ is $\an$. Equivalently, $S_\vp$ is $\an$ if and only if $T_{\bar{\vp}}$ is $\an$,
\item[(ii)] $T_\vp$ is $\am$ if and only if $S_{\bar{\vp}}$ is $\am$,
\item[(iii)] $T_\vp \in \overline{\an(H^2)}$ if and only if $S_{\bar{\vp}}\in \overline{\an(H^2_{-})}$. 
\end{enumerate}
\end{remark}

The space of quasi continuous functions is defined by $\mathcal{QC}=[H^\infty + C(\T)] \cap [\overline{H^\infty + C(\T)}]$. It is a uniformly closed self-adjoint sub-algebra of $L^\infty(\T)$. See \cite[Page 157]{Douglas}, \cite[Page 729]{Peller} for further details.

\begin{example}
    Let $\vp \in L^\infty (\T)$ and $T=\begin{bmatrix}
        T_\vp & 0 \\
        0 & S_\vp
    \end{bmatrix}$ on $H^2\oplus H^2_{-}$. From \cite[Theorem 8]{CR Paris}, we have $T_\vp, T_{\bar{\vp}} \in \an (H^2)$ if and only if $|\vp|=\|\vp\|_{\infty}$ a.e. $(\mu)$ and $(I-P)\vp$, $(I-P)\bar{\vp}$ are rational functions. Now $T^*T$ is $\an$ if and only if $T^*_\vp T_\vp, S^*_\vp S_\vp$ are $\an$ with $\es(T^*_\vp T_\vp)=\es(S^*_\vp S_\vp)$, by Lemma \ref{Diagonal AN}. Thus, by Lemma \ref{AN T^*T} and Remark \ref{Anti}, we get $T\in \an(H^2\oplus H^2_{-})$ if and only if $|\vp|=\|\vp\|_{\infty}$ a.e. and $(I-P)\vp,(I-P)\bar{\vp}$ are rational functions. Note that applying Weyl's theorem in  \cite[Theorem 2.1]{AN Toeplitz}, $\es(|T_\vp|)=\es(|T_{\bar{\vp}}|) = \{\|\vp \|_{\infty} \}$ and by Remark \ref{Ess TO and DTO}, $\es(|T_\vp|)=\es(|S_\vp|)$. Thus, the condition $\es(T^*_\vp T_\vp)=\es(S^*_\vp S_\vp)$ follows automatically.
    
    Now $T_\vp, T_{\bar{\vp}} \in \am (H^2)$ if and only if $|\vp|=\alpha$ a.e., where $\alpha \in \es(|T_\vp|)=\es(|T_{\bar{\vp}}|)$ and $\vp,\bar{\vp} \in H^\infty + C(\T)$, by \cite[Theorem 12]{CR Paris}. Note that $\vp, \bar{\vp}\in H^\infty + C(\T)$ if and only if $\vp \in \mathcal{QC}$. Proceeding with the similar steps as in $\an$ case, we can show that $T\in \am(H^2\oplus H^2_{-})$ if and only if $\vp \in \mathcal{QC}$ with $|\vp|=\alpha$ a.e., where $\alpha \in \es(|T_\vp|)=\es(|S_\vp|)$.\\
    By Theorem \ref{AN Closure Matrix}, $T^*T=\begin{bmatrix}
        T^*_\vp T_\vp & 0 \\
        0 & S^*_\vp S_\vp
    \end{bmatrix} \in \overline{\an (L^2(\T))}_{+}$ implies $T^*_\vp T_\vp \in \overline{\an (H^2)}_{+}$, $S^*_\vp S_\vp \in \overline{\an (H^2_{-})}_{+}$ and $\es(T^*_\vp T_\vp)=\es(S^*_\vp S_\vp)$. Hence by Lemma \ref{AN Closure T^*T} and Remark \ref{Anti}, $T_\vp , T_{\bar{\vp}} \in \overline{\an (H^2)}$. From  \cite[Corollary 2.4]{AN Toeplitz}, we can say $|\vp|=\alpha$, constant a.e. $(\mu)$ and further $\es(T^*_\vp T_\vp)=\{\alpha^2\}=\es(S^*_\vp S_\vp)$. Thus, by Lemma \ref{AN Closure T^*T}, if $T\in \overline{\an(L^2(\T))}$, then $|\vp|=\alpha$, constant a.e. $(\mu)$ and $\es(|T_\vp|)=\es(|S_\vp|)=\{\alpha\}$.
\end{example}

\begin{example}
    Suppose $\vp \in L^\infty (\T)$ and $T=\begin{bmatrix}
        I_{H^2} & H^*_\vp \\
        H_\vp & I_{H^2_{-}}
    \end{bmatrix}$ on $H^2 \oplus H^2_{-}$. Here $T$ is self-adjoint and $T^2=\begin{bmatrix}
        I_{H^2} + H^*_\vp H_\vp & 2 H^*_\vp \\
        2 H_\vp  &  I_{H^2_{-}} + H_\vp H^*_\vp
    \end{bmatrix}$. If $T^2\in \an (L^2(\T))_{+}$, then by Theorem \ref{AN Necessary}, $H^*_\vp$ is compact i.e. $H_\vp$ is compact. By \cite[Theorem 5.5, Page 27]{Peller}, we get $\vp \in H^\infty + C(\T)$. If $(I-P)\vp$ is rational function, $H_\vp$ is of finite rank, by \cite[Corollary 3.2, Page 21]{Peller}. Then clearly the diagonal entries of $T^2$ are $\an$ with equal essential spectrum. Thus, by Theorem \ref{AN Sufficient}, $T^2$ is $\an$ and hence $T$ is $\an$.
    
    By Theorem \ref{AN Closure Matrix}, $T^2\in \overline{\an (L^2(\T))}$ if and only if $H^*_\vp$ is compact i.e. $H_\vp$ is compact. Here the compactness of the operator $H_\vp$ is enough as the remaining conditions of Theorem \ref{AN Closure Matrix} follows trivially. Hence $T \in \overline{\an (L^2(\T))}$ if and only if $\vp \in H^\infty + C(\T)$.

   If  $\|\vp \| \leq 1$, $\|H_\vp \|=dist (\vp,H^\infty) \leq \|\vp \|\leq 1$ i.e. $H_\vp$ is a contraction. By Lemma \ref{AN Contraction}, $T\in \an (H^2 \oplus H^2_{-})$ if and only if $H^*_\vp \in \clf (H^2_{-},H^2)$, if and only if $H_\vp \in \clf (H^2,H^2_{-})$ if and only if $(I-P)\vp$ is a rational function.
\end{example}

\begin{example}
    Let $\vp \in L^\infty (\T)$ and consider the multiplication operator $M_\vp$ on $L^2(\T)$. We have $M_\vp=\begin{bmatrix}
        T_\vp & H^*_{\bar{\vp}} \\
        H_\vp & S_\vp
    \end{bmatrix}$ on $H^2 \oplus H^2_{-}=L^2(\T)$. If $M_\vp$ is $\an$, then $M^*_\vp M_\vp$ is also $\an$. Now $M^*_\vp M_\vp =M_{|\vp|^2}=\begin{bmatrix}
        T_{|\vp|^2} & H^*_{|\vp|^2} \\
        H_{|\vp|^2} & S_{|\vp|^2}
    \end{bmatrix}$. From Theorem \ref{AN Necessary}, we get $T_{|\vp|^2}, S_{|\vp|^2}$ are $\an$ with $\es (T_{|\vp|^2})=\es(S_{|\vp|^2})$. From the relation $T_{|\vp|^2}=VS_{|\vp|^2}V$, $\an$ condition for $S_{|\vp|^2}$ is same as $T_{|\vp|^2}$ has. By \cite[Theorem 8]{CR Paris}, we get $|\vp|^2=\||\vp|^2 \|_{\infty}$ a.e. $(\mu)$, which further implies $|\vp|=\|\vp \|_{\infty}$ a.e. $(\mu)$. Here the off-diagonal entries of $M^*_\vp M_\vp$ become zero as $|\vp|$ is inner. Also, $\es(M^*_\vp M_\vp )=\es(T_{|\vp|^2})=\es(S_{|\vp|^2})=\{\|\vp\|^2_{\infty}\}$.
\end{example}

\begin{example}
Let $\vp \in L^\infty(\T)$ be a real function. Then $M_\vp$ is a self-adjoint operator. Also, $M^2_\vp = M_{\vp^2}=\begin{bmatrix}
    T_{\vp^2} & H^*_{\vp^2} \\
    H_{\vp^2} & S_{\vp^2}
\end{bmatrix}$ is positive. If $M_{\vp^2}$ is $\an$, then by Theorem \ref{AN Necessary}, $T_{\vp^2}$ and $S_{\vp^2}$ are $\an$. From Remark \ref{Anti}, the $\an$ condition for $S_{\vp^2}$ is exempted by $T_{\vp^2}$. Thus, $\vp^2=\|\vp^2\|_{\infty}$ a.e. $(\mu)$, by \cite[Theorem 8]{CR Paris}. Consequently, we get $\vp =\pm \|\vp\|$ a.e. $(\mu)$. Similarly if $M_{\vp^2}$ is $\am$, then using Theorem \ref{AM Necessary}, Remark \ref{Anti} and \cite[Theorem 12]{CR Paris}, we get $\vp^2=\alpha$, constant a.e. $(\mu)$, that is $\vp=\pm \sqrt{\alpha}$ a.e. $(\mu)$. By Lemma \ref{AN T^*T} and Lemma \ref{AM T^*T}, we conclude that if $M_\vp$ is $\an$ or $\am$, then $\vp$ is constant almost everywhere on $\T$.
\end{example}

\begin{remark}
    If $\vp\in L^\infty(\T)$ is a non-constant real function, then $M_\vp$ is neither $\an$ nor $\am$.
\end{remark}

\begin{center}
		\textbf{Statements and Declarations}
	\end{center}
	
		\noindent \textbf{Data Availability}\\
	Data sharing is not applicable to this article, as no data sets were generated or
analyzed during the current study.\\
	
	\noindent \textbf{Conflict of interest}\\
	On behalf of all authors, the corresponding author states that there is no conflict of interest. \\

\end{document}